\documentclass[12pt]{article}
\usepackage{amsmath,amssymb,amsthm}

\usepackage{hyperref}

\usepackage{multicol}

\newtheorem{theorem}{Theorem}[section]
\newtheorem{lemma}{Lemma}[section]
\newtheorem{corollary}{Corollary}[section]
\newtheorem{definition}{Definition}[section]

\begin{document}

\title{A Density Bound for Mu{\~n}iz's List Reversal}

\author{Alan Talmage\footnote{atalmage@semo.edu \newline Southeast Missouri State University \newline Cape Girardeau, Missouri, 63701, USA}}

\maketitle

\begin{abstract}
Mu{\~n}iz described the problem of reversing a list of positive integers under the following rules: given a list of distinct positive integers, at each step, one can combine two adjacent entries into their sum or break an entry into two entries with the same sum, under the constraints that after each step the list must contain distinct entries no larger than the original largest entry.  The question of what lists are and are not reversible under these rules is studied and an asymptotic upper bound on the density of an irreversible list is established.
\end{abstract}


\section{Introduction}
In 2024, Mu{\~n}iz introduced a novel combinatorial puzzle~\cite{muniz}, which he called Reverse the List of Numbers.  While this puzzle sparked a flurry of interest at the time (see, e.g., discussions at~\cite{muniz}, ~\cite{hn}, implementations and solvers at~\cite{lawson-perfect},~\cite{alayachew} ~\cite{61u}), much of the asymptotic behavior of the problem's solvability has remained unexamined.  The problem statement is as follows:

Given an ordered list $L$ of distinct positive integers with maximal entry $n$, we allow the following operations:
\begin{itemize}
\item \emph{Split:} A list entry may be split into two entries whose sum is the original entry,
\item \emph{Combine:} Two adjacent entries in the list may be combined into a single entry equal to their sum,
\end{itemize}
where the list must satisfy the following conditions after each step:
\begin{itemize}
\item The list entries are distinct,
\item The largest entry in the list is no larger than $n$, the original maximum entry.
\end{itemize}

The goal is to take a list and transform it into its reverse using only these operations.  We call a list \emph{reversible} if that is possible, and \emph{irreversible} if not.

For example, the following is a valid reversal of the list $(3, 5, 7)$ (as seen in~\cite{muniz}).

\begin{center}
$3, \hspace{3ex} 5, \hspace{3ex} 7$

$2, \hspace{3ex} 1, \hspace{3ex} 5, \hspace{3ex} 7$

$2, \hspace{3ex} 1, \hspace{3ex} 5, \hspace{3ex} 4, \hspace{3ex} 3$

$2, \hspace{3ex} 6, \hspace{3ex} 4, \hspace{3ex} 3$

$2, \hspace{3ex} 5, \hspace{3ex} 1, \hspace{3ex} 4, \hspace{3ex} 3$

$7, \hspace{3ex} 1, \hspace{3ex} 4, \hspace{3ex} 3$

$7, \hspace{3ex} 5, \hspace{3ex} 3$
\end{center}

Our goal is to examine conditions conditions guaranteeing reversibility.  Intuitively, a list of very low density should be reversible, due to the large number of possible moves available.  Following this intuition, we prove an explicit lower bound on the density of an irreversible list.  

To approach the problem, we begin by studying the ability to reorder a list at will, which is a stronger condition than irreversibility.  The main result is the following:

\begin{theorem}\label{main_theorem}
Let $L$ be a list of $\ell$ entries with maximal entry $n$.  If $n > 5\ell +3 $, any two adjacent list entries can be transposed.
\end{theorem}

That is, there exists a valid sequence of operations from the list $(\dots, a, b, \dots)$ to the list $(\dots, b, a, \dots)$, with the rest of the list returned to its original state.

Since the ability to transpose arbitrary adjacent list entries allows us to reorder the list at will, we have an immediate corollary.

\begin{corollary}\label{main_corollary}
Let $L$ be a list of $\ell$ entries with maximal entry $n$.  If $n > 5\ell + 3$, the list is reversible.
\end{corollary}

Given a fixed length $\ell$, there are only finitely many lists with $n \le 5\ell + 3$, so we can observe:

\begin{corollary}\label{finite_lists}
There are finitely many irreversible lists of each length.
\end{corollary}

A \emph{reversal} will be a sequence of valid operations that transform a list into it its reverse.  A \emph{swap} will be a sequence of valid operations that transpose two adjacent list entries.  We call a number \emph{open} or \emph{unobstructed} if it is does not currently appear in the list and \emph{obstructed} if it does.  To avoid counting a list and its reverse twice, we will assume that the first entry of a list is smaller than the last entry when counting lists.  We assume throughout this work that a list $L$ has length $\ell$ and largest element $n$.

\section{Preliminaries: Lists of Length 2}

We begin by examining lists $(a, b)$ of length 2.  There are two general tactics here: if $a$ and $b$ are widely spaced, then we can split value from $b$ and combine it into $a$ in two steps to swap them,
%
%
%
%
while if $a$ and $b$ are very close in size, we can split them both into smaller numbers and combine these together into the reversed order.  These two general approaches to swapping a pair of entries will be reoccurring techniques throughout this work.

\begin{lemma}\label{length_two_lists}
There are eleven irreversible length 2 lists, which are 
\[
	(1,2), (1,3), (1,4), (1,5), (2,3), (2,4), (2,5), (2,6), (3,4), (3,5), (4,5).
\]
\end{lemma}
\begin{proof}
The irreversibility of these eleven lists can be verified directly, so we focus one showing that all other lists of length 2 are reversible.  

Consider the list $(a, b)$.  The swap 
\begin{center}
$a, \hspace{3ex} b$

$a, \hspace{3ex} k, \hspace{3ex} b-k$

$a+k, \hspace{3ex} b+k$

$a+k, \hspace{3ex} b-k-a, \hspace{3ex} a$

$b, \hspace{3ex} a$
\end{center}

will reverse all lists except for a few with specific spacings between $a$ and $b$.  We begin by examining the $k=1$ and $k=2$ cases of this swap.

For $k=1$,
\begin{center}
$a, \hspace{3ex} b$

$a, \hspace{3ex} 1, \hspace{3ex} b-1$

$a+1, \hspace{3ex} b+1$

$a+1, \hspace{3ex} b-1-a, \hspace{3ex} a$

$b, \hspace{3ex} a$
\end{center}

reverses the list $(a, b)$ if $1 < a < b$, $2a \ne b-1$, and $a+1 \ne b+1$.  Thus the only cases remaining are lists of the form $(1, b), (a, a+1), (a, a+2), (a, 2a+1), (a, 2a+2)$.

For $k=2$,
\begin{center}
$a, \hspace{3ex} b$

$a, \hspace{3ex} 2, \hspace{3ex} b-2$

$a+2, \hspace{3ex} b-2$

$a+2, \hspace{3ex} b-2-a, \hspace{3ex} a$

$b, \hspace{3ex} a$
\end{center}

reverses the list $(a, b)$ if $a \ne 2$, $b > 2+a$, $b \ne a+4$, $b \ne 2a+2$.  This handles all cases except lists of the form $(2, b), (a, a+1), (a, a+2), (a, a+4), (a, 2a+2)$.

Comparing the two, we are left with the cases $(a, a+1), (a, a+2), (a, 2a+2)$ that cannot be reversed by either of the above approaches.  We can now examine swaps for each of these remaining cases.

The swap
\begin{center}
$2, \hspace{3ex} b$

$2, \hspace{3ex} b-3, \hspace{3ex} 3$

$b-1, \hspace{3ex} 3$

$b-1, \hspace{3ex} 1, \hspace{3ex} 2$

$b, \hspace{3ex} 2$
\end{center}

reverses $(2, b)$ for $b \ge 7$,

\begin{center}
$a, \hspace{3ex} a+1$

$2, \hspace{3ex} a-2, \hspace{3ex} a+1$

$2, \hspace{3ex} a-2, \hspace{3ex} 1, \hspace{3ex} a$

$2, \hspace{3ex} a-1, \hspace{3ex} a$

$a+1, \hspace{3ex} a$
\end{center}

reverses $(a, a+1)$ for $a \ge 5$,

\begin{center}
$a, \hspace{3ex} a+2$

$1, \hspace{3ex} a-1, \hspace{3ex} a+2$

$1, \hspace{3ex} a-1, \hspace{3ex} 2, \hspace{3ex} a$

$1, \hspace{3ex} a+1, \hspace{3ex} a$

$a+2, \hspace{3ex} a$
\end{center}

reverses $(a, a+2)$ for $a \ge 4$, and

\begin{center}
$a, \hspace{3ex} 2a+2$

$a, \hspace{3ex} 3, \hspace{3ex} 2a-1$

$a+3, \hspace{3ex} 2a-1$

$a+3, \hspace{3ex} a-1, \hspace{3ex} a$

$2a+2, \hspace{3ex} a$
\end{center}

reverses $(a, 2a+2)$ for $a \ge 5$.

Directly verifying the finite number of lists too small for these swaps yields the lemma.

\end{proof}

Note that all of these irreversible lists of length 2 all violate the $n > 5\ell + 3$ assumption of Theorem~\ref{main_theorem}, so this lemma implies the theorem in the case of a list of length 2.

\section{A Density Bound}

We are now ready for the proof of our main theorem.

\begin{proof}[Proof of Theorem~\ref{main_theorem}]
Let $L$ be a list of $\ell$ entries with maximal entry $n$ satisfying $n > 5\ell + 3$.  Let $L = (\dots, a, b, \dots)$, with the understanding that either ellipsis may be empty, and recall that our goal is to reach the list $(\dots, b, a, \dots)$, which is identical to $L$ except for the transposition of $a$ and $b$.  We may assume without loss of generality that $b > a$.

We divide the proof into several cases based on $a$ and $b$.

\begin{itemize}
\item Case 1: $a + b \le n$
\item Case 2: $a < \frac{n}{2}$
\begin{itemize}
\item Case 2a: $2a+b < 2n - 2\ell - 1$
\item Case 2b: $a < \frac{n}{2}$, $b > 2\ell$
\end{itemize}
\item Case 3: $a \ge \frac{n}{2}$
\end{itemize}

These cases will be handled in order, with each case dependent on the previous cases.  Note that Cases 2a and 2b together cover Case 2, since if $a < \frac{n}{2}$ and $b \le 2\ell$, then $2a + b < 2\frac{n}{2} + 2\ell = n + 2\ell < n + (5\ell + 3) - 2\ell - 1 < 2n - 2\ell - 1$, so any instance with $a < \frac{n}{2}$ is contained in either Case 2a or Case 2b.

\subsection*{Case 1: $a + b \le n$}

The goal for this case is to use the swap

\begin{center}
$a, \quad b$

$a+b$

$b, \quad a$
\end{center}

to reverse the pair.  If $a+b$ is open, this immediately works.  If it is obstructed, we make it the first element what we will call a chain of obstructions.
\begin{definition}
A \emph{chain of obstructions} for a list $L$ is an ordered sequence $\{x_1, x_2, \dots\}$ such that
\begin{itemize}
\item each $x_i$ is obstructed (i.e. is an entry in $L$), and 
\item if the value $x_i$ were made open, then we can open the value $x_{i-1}$.
\end{itemize}
\end{definition}

Consider the element $c$ that is currently last in the chain of obstructions (beginning with $c = a+b$).  If $c > 2\ell$, then there are $> \left\lfloor\frac{c}{2}\right\rfloor \ge \ell$ ways to split it.  However, $c$ itself cannot obstruct any of these, so there are at most $\ell - 1$ obstructions for $\ell$ possible splits, so we can split it, making the value $c$ open and thus collapsing the chain.  

Otherwise, consider an element $c$ of the chain satisfying $c \le 2\ell$ that has a neighbor $d$ that is not already in the chain.  If $c+d \le n$, then we can extend the chain by appending $c+d$ to it.  Otherwise, we may assume that $d \ge n - c > n - 2\ell$.

The manipulation

\begin{center}
$c, \hspace{3ex} d$

$c, \hspace{3ex} k, \hspace{3ex} d-k$

$c+k, \hspace{3ex} d-k$
\end{center}

will open the value $c$.  This requires $k < d$, with no two entries in the same step equal:
\[
	\quad k \ne c,  c \ne d-k, k \ne d-k, c+k \ne d-k, 
\]
with the three values $k, d-k, c+k$ unobstructed.  There are $> n - 2\ell$ possible choices of $k$.  At most four are removed by the conditions avoiding duplicate entries in each step, and each of the $\ell - 2$ list entries that are not $c$ or $d$ can obstruct at most three values of $k$, so there are
\[
	\ge n - 2\ell - 4 - 3(\ell-2) = n - 5\ell + 2
\]
values of $k$ that open $c$.  Such a $k$ exists if
\[
	n - 5\ell + 2 > 0.
\]

Thus if $n > 5\ell - 2$, the chain of obstructions must eventually collapse, which makes $a+b$ open, at which point we can swap $a, b$.

\subsection*{Case 2a: $2a+b < 2n - 2\ell - 1$}

The idea here is to split $b$ into $k, b-k$ such that $a + k \le n$ and $a + b - k \le n$.  In this case, we can use Case 1 (with $\ell$ replaced by $\ell + 1$) to swap $a, k$ and $a, b-k$, allowing us to perform the swap

\begin{center}
$a, \hspace{3ex} b$

$a, \hspace{3ex} k, \hspace{3ex} b-k$

$k, \hspace{3ex} a, \hspace{3ex} b-k$

$k, \hspace{3ex} b-k, \hspace{3ex} a$

$b, \hspace{3ex} a$.
\end{center}

This requires a $k$ with $k < b$, $k \ne a, k \ne b - a, k \ne \frac{b}{2}$, with $k$ and $b-k$ open, $a + k \le n$, $a + b - k \le n$, and Case 1 to accept $\ell + 1$ in place of $\ell$.

Thus we have $a + b - n \le k < \min(n-a, b)$ with three forbidden values of $k$ and two numbers that must be unobstructed.  Since $a + b \le n$ has been handled by Case 1, we may assume $a + b > n$, so $b > n - a$, so $\min(n-a, b) = n - a$.  Thus the number of $k$ satisfying this is
\[
	\ge n - a - (a + b - n) - 1 - 2\ell.
\]

We need this to be positive, so we require
\[
	n - a - (a + b - n) - 1 - 2\ell > 0
\]
\[
	\Leftrightarrow 2a + b < 2n - 2\ell - 1.
\]

In addition, Case 1 required $n > 5\ell - 2$ and we need it to work on a list of length $\ell + 1$, so we also require
\[
	n > 5(\ell + 1) - 2 = 5\ell + 3.
\]

\subsection*{Case 2b: $a < \frac{n}{2}$, $b > 2\ell$}

The approach here is to split $b$ into two parts $k, b - k$ such that Case 2a applies to both $(a, k)$ and $(a, b-k)$.  Then $a$ and $b$ can be swapped by the process

\begin{center}
$a, \hspace{3ex} b$

$a, \hspace{3ex} k, \hspace{3ex} b-k$

$k, \hspace{3ex} a, \hspace{3ex} b-k$

$k, \hspace{3ex} b-k, \hspace{3ex} a$

$b, \hspace{3ex} a$
\end{center}

where the second and third operations are done by Case 2a.

This requires $(a, k)$ and $(a, b-k)$ to meet the conditions of Case 2a with $\ell$ replaced by $\ell+1$.  Since interchanging the roles of $a$ and $b$ in the bound $2a+b < 2n - 2\ell - 1$ for Case 2a only makes the bound stronger, we can do so at will here, so it suffices to have $2k + a < 2n - 2(\ell+1) - 1 = 2n - 2\ell - 3$ and $2k + (b-k) < 2n - 2\ell - 3$.  These are equivalent to $k < n - \ell - \frac{3}{2} - \frac{a}{2} < \frac{3n}{4} - \ell - \frac{3}{2}$ and $k < 2n - 2\ell - 3 - b < n - 2\ell - 3$.

As long as $n > 4\ell$, $\frac{3n}{4} < n - \ell$, so both of these inequalities are satisfied if we have $k < \frac{3n}{4} - \ell - 3$.  Since we are splitting $b$ into $k$ and $b-k$, we can assume without loss of generality that $k < \frac{b}{2}$, and as long as $n > 4\ell + 3$, we have $\frac{b}{2} < \frac{n}{2} < \frac{3n}{4} - \ell - 3$, so $k < \frac{b}{2}$ suffices.

We require $k \ne a$ and $b - k \ne a$ with $k$, $b-k$ unobstructed, i.e. $k \ne b - a$ with $k, b-k$ unobstructed by any of the $\le \ell-1$ elements of $L$ that are $<b$.  Thus the number of such $k$ is
\[
	\ge \left(\frac{b}{2} - 1\right) - 1 - (\ell-1) = \frac{b}{2} - \ell.
\]

This is $> 0$ if $b > 2\ell$, which concludes this case.

\subsection*{Case 3: $a \ge \frac{n}{2}$}

Cases 1, 2a, and 2b now allow us to swap $a, b$ as long as $a < \frac{n}{2}$.  As long as $a$ is not too much larger, we should be able to split it into two parts and apply the previous cases.  We can then apply this technique recursively.

Let $A_i(\ell)$ be a value such that if $a < A_i(\ell)$, then $a, b$ can be swapped by splitting $a$ recursively $i$ times, with $A_0(\ell) = \frac{n}{2}$ being that which is handled by Cases 1, 2a, and 2b.  Then, to find $A_{i+1}(\ell)$, we consider splitting

\begin{center}
$a$

$k, \hspace{3ex} a-k$
\end{center}

with $k, a-k < A_i(\ell)$ requires $k \ne a-k$ with $a, a-k$ unobstructed.  If we assume that $k$ is the smaller part, i.e., $k < \frac{a}{2}$, this is equivalent to $a - A_i(\ell) < k < \frac{a}{2}$ with $k, a-k$ unobstructed.  There are $\ge \lfloor\frac{a}{2}\rfloor - (a - A_i(\ell)) - 1 \ge\frac{a}{2} - a + A_i(\ell) - \frac{3}{2}$ such $k$, with at most $\ell - 2$ elements of $L$ smaller than $a$ (thus not $a$ or $b$) obstructions.  The number of valid $k$ is then
\[
	\ge \frac{a}{2} - a + A_i(\ell) - \frac{3}{2} - (\ell-2)
\]
\[
	= A_i(\ell) - \frac{a}{2} - \ell + \frac{1}{2}.
\]

This is $>0$ when 
\[
	a < 2A_i(\ell) - 2\ell + 1,
\]

Now we lengthen the list by one when we split $a$, which means that this gives the recursion
\[
	A_{i+1}(\ell) = 2A_i(\ell+1) - 2\ell + 1.
\]

With the initial condition $A_0(\ell)$, we can calculate

\begin{center}
\begin{tabular}{|c|c|}
$i$ & $A_i(\ell)$ \\
\hline
0 & $\frac{n}{2}$ \\
1 & $n - 2\ell + 1$ \\
2 & $2n - 6\ell - 1$ \\
3 & $4n - 14\ell - 13$ \\
\end{tabular}
\end{center}

We need to find $A_i(\ell) > n$, since then any value of $b$ is included.  At $i=3$ we have 
\[
	A_3(\ell) = 4n - 14\ell - 13 > n
\]
\[
	\Leftrightarrow n > \frac{14}{3} \ell + \frac{13}{3}.
\]

Since $n > 5\ell + 3$ is an integer, we have $n > 5\ell + \frac{11}{3}$, so
\[
	n > 5\ell + \frac{11}{3} = \frac{14}{3}\ell + \frac{\ell}{3} + \frac{11}{3} \ge \frac{14}{3}\ell + 1 + \frac{11}{3} > \frac{14}{3} \ell + \frac{13}{3}, 
\]

which completes this case, finishing the proof.

\end{proof}

%
%
%
%
%
%
%
%
%

\section{Discussion}
The $n > 5\ell + 3$ bound of Theorem~\ref{main_theorem} is almost certainly not the best possible density bound.  The lowest-density irreversible list known to the author is $[2, 6]$ with density $\frac{n}{\ell} = \frac{1}{3}$.  It appears probable that that is the best such example, and that the best possible bound would be $n > 3\ell$.  In fact, the lowest-density irreversible list for $\ell \ge 3$ appears to be $\frac{1}{2}$.

%


There are several other interesting problems raised by the Reverse the List of Numbers puzzle.  A few such questions are listed here to spark the imagination of the curious reader.

\begin{itemize}

\item What is the computational complexity of the decision problem for list reversibility?  (A graph search can answer the decision problem, but it is not immediately apparent how large the graph may be relative to the input length.)

\item Define the \emph{swap time} of a list to be the minumum number of swaps needed to reverse it.  What is the maximum swap time among all lists of length $\ell$, and how does this maximum increase with $\ell$?  (See~\cite{rokicki} for some data on early values of the maximum swap time, including the fact that it is not monotonic in $\ell$.)

\item The irreversible lists found by computer search in~\ref{irreversible_lists} suggest a conjecture that reversibility is monotonic in $n$, in the sense that replacing the maximum value $n$ in an irreversible list with $n+1$ (leaving the rest of the list unchanged) must yield an irreversible list.  No counterexamples to the conjecture are presently known, but it is not clear if that is an artifact of only having data on short lists.

\item One could also study the question of reversibility where the requirement that all list entries are smaller than a fixed $n$ at each stage does not require $n$ to be the maximum value of the original list.  (This situation naturally arises at intermediate steps of a list reversal.)  Notably, this form of reversibility is not monotonic in the sense given above.  For example, the list $(1,2,7,3)$ is reversible under maximum value 9, but the list $(1,2,8,3)$ is not.  

\end{itemize}

\appendix
\section{Known Irreversible Lists}\label{irreversible_lists}

The following are irreversible lists found by brute-force computer search.  For each length, lists with maximum value up to the listed $n$ have been searched.  It has not yet been proved that this listing contains all irreversible lists of these lengths, but that is conjectured to be the case.

\subsection*{$\ell = 3$, $n \le 15$: 51 known irreversible lists}

\begin{multicols}{4}

[1, 2, 3]

[1, 2, 4]

[1, 2, 5]

[1, 3, 2]

[1, 3, 4]

[1, 3, 5]

[1, 4, 2]

[1, 4, 3]

[1, 4, 5]

[1, 5, 2]

[1, 5, 3]

[1, 5, 4]

[1, 6, 2]

[1, 6, 4]

[1, 6, 5]

[1, 7, 5]

[2, 1, 3]

[2, 1, 4]

[2, 1, 5]

[2, 3, 4]

[2, 3, 5]

[2, 4, 3]

[2, 4, 5]

[2, 5, 3]

[2, 5, 4]

[2, 5, 6]

[2, 6, 3]

[2, 6, 4]

[2, 6, 5]

[2, 7, 4]

[3, 1, 4]

[3, 1, 5]

[3, 2, 4]

[3, 2, 5]

[3, 4, 5]

[3, 4, 6]

[3, 5, 4]

[3, 5, 6]

[3, 6, 4]

[3, 6, 5]

[4, 1, 5]

[4, 2, 5]

[4, 3, 5]

[4, 3, 6]

[4, 5, 6]

[4, 6, 5]

[4, 7, 5]

[4, 7, 6]

[5, 2, 6]

[5, 3, 6]

[5, 4, 6]

\end{multicols}

\subsection*{$\ell = 4$, $n \le 10$: 307 known irreversible lists}

\begin{multicols}{4}

[1, 2, 3, 4]

[1, 2, 3, 5]

[1, 2, 4, 3]

[1, 2, 4, 5]

[1, 2, 4, 6]

[1, 2, 5, 3]

[1, 2, 5, 4]

[1, 2, 5, 6]

[1, 2, 6, 3]

[1, 2, 6, 4]

[1, 2, 6, 5]

[1, 2, 7, 3]

[1, 2, 8, 3]

[1, 3, 2, 4]

[1, 3, 2, 5]

[1, 3, 4, 2]

[1, 3, 4, 5]

[1, 3, 4, 6]

[1, 3, 5, 2]

[1, 3, 5, 4]

[1, 3, 5, 6]

[1, 3, 6, 2]

[1, 3, 6, 4]

[1, 3, 6, 5]

[1, 3, 7, 2]

[1, 3, 7, 5]

[1, 3, 7, 6]

[1, 3, 8, 4]

[1, 4, 2, 3]

[1, 4, 2, 5]

[1, 4, 2, 6]

[1, 4, 3, 2]

[1, 4, 3, 5]

[1, 4, 3, 6]

[1, 4, 5, 2]

[1, 4, 5, 3]

[1, 4, 5, 6]

[1, 4, 6, 2]

[1, 4, 6, 3]

[1, 4, 6, 5]

[1, 4, 7, 5]

[1, 5, 2, 3]

[1, 5, 2, 4]

[1, 5, 2, 6]

[1, 5, 3, 2]

[1, 5, 3, 4]

[1, 5, 3, 6]

[1, 5, 4, 2]

[1, 5, 4, 3]

[1, 5, 4, 6]

[1, 5, 6, 2]

[1, 5, 6, 3]

[1, 5, 6, 4]

[1, 5, 7, 4]

[1, 6, 2, 3]

[1, 6, 2, 4]

[1, 6, 2, 5]

[1, 6, 3, 2]

[1, 6, 3, 4]

[1, 6, 3, 5]

[1, 6, 4, 2]

[1, 6, 4, 3]

[1, 6, 4, 5]

[1, 6, 5, 2]

[1, 6, 5, 3]

[1, 6, 5, 4]

[1, 7, 2, 3]

[1, 7, 3, 2]

[1, 7, 3, 5]

[1, 7, 3, 6]

[1, 7, 4, 5]

[1, 7, 4, 6]

[1, 7, 5, 3]

[1, 7, 5, 4]

[1, 7, 5, 6]

[1, 7, 6, 2]

[1, 7, 6, 3]

[1, 7, 6, 4]

[1, 7, 6, 5]

[1, 8, 6, 4]

[1, 8, 7, 3]

[2, 1, 3, 4]

[2, 1, 3, 5]

[2, 1, 4, 3]

[2, 1, 4, 5]

[2, 1, 4, 6]

[2, 1, 5, 3]

[2, 1, 5, 4]

[2, 1, 5, 6]

[2, 1, 6, 3]

[2, 1, 6, 4]

[2, 1, 6, 5]

[2, 1, 7, 3]

[2, 1, 8, 3]

[2, 3, 1, 4]

[2, 3, 1, 5]

[2, 3, 4, 5]

[2, 3, 4, 6]

[2, 3, 5, 4]

[2, 3, 5, 6]

[2, 3, 6, 4]

[2, 3, 6, 5]

[2, 3, 7, 4]

[2, 3, 7, 5]

[2, 4, 1, 3]

[2, 4, 1, 5]

[2, 4, 1, 6]

[2, 4, 3, 5]

[2, 4, 3, 6]

[2, 4, 5, 3]

[2, 4, 5, 6]

[2, 4, 6, 3]

[2, 4, 6, 5]

[2, 5, 1, 3]

[2, 5, 1, 4]

[2, 5, 1, 6]

[2, 5, 3, 4]

[2, 5, 3, 6]

[2, 5, 4, 3]

[2, 5, 4, 6]

[2, 5, 6, 3]

[2, 5, 6, 4]

[2, 5, 7, 6]

[2, 6, 1, 3]

[2, 6, 1, 4]

[2, 6, 1, 5]

[2, 6, 3, 4]

[2, 6, 3, 5]

[2, 6, 4, 3]

[2, 6, 4, 5]

[2, 6, 4, 7]

[2, 6, 5, 3]

[2, 6, 5, 4]

[2, 6, 5, 7]

[2, 6, 7, 3]

[2, 6, 7, 4]

[2, 6, 7, 5]

[2, 7, 1, 3]

[2, 7, 3, 5]

[2, 7, 3, 6]

[2, 7, 4, 5]

[2, 7, 4, 6]

[2, 7, 5, 3]

[2, 7, 5, 4]

[2, 7, 5, 6]

[2, 7, 6, 3]

[2, 7, 6, 4]

[2, 7, 6, 5]

[2, 7, 8, 3]

[2, 7, 8, 6]

[2, 8, 3, 6]

[2, 8, 4, 5]

[2, 8, 5, 4]

[2, 8, 6, 3]

[2, 8, 6, 7]

[2, 8, 7, 6]

[3, 1, 2, 4]

[3, 1, 2, 5]

[3, 1, 4, 5]

[3, 1, 4, 6]

[3, 1, 5, 4]

[3, 1, 5, 6]

[3, 1, 6, 4]

[3, 1, 6, 5]

[3, 1, 7, 5]

[3, 1, 7, 6]

[3, 2, 1, 4]

[3, 2, 1, 5]

[3, 2, 4, 5]

[3, 2, 4, 6]

[3, 2, 5, 4]

[3, 2, 5, 6]

[3, 2, 6, 4]

[3, 2, 6, 5]

[3, 2, 7, 4]

[3, 2, 7, 5]

[3, 4, 1, 5]

[3, 4, 1, 6]

[3, 4, 2, 5]

[3, 4, 2, 6]

[3, 4, 5, 6]

[3, 4, 6, 5]

[3, 4, 7, 6]

[3, 5, 1, 4]

[3, 5, 1, 6]

[3, 5, 2, 4]

[3, 5, 2, 6]

[3, 5, 4, 6]

[3, 5, 4, 7]

[3, 5, 6, 4]

[3, 5, 6, 7]

[3, 5, 7, 4]

[3, 5, 7, 6]

[3, 6, 1, 4]

[3, 6, 1, 5]

[3, 6, 2, 4]

[3, 6, 2, 5]

[3, 6, 4, 5]

[3, 6, 4, 7]

[3, 6, 5, 4]

[3, 6, 5, 7]

[3, 6, 7, 4]

[3, 6, 7, 5]

[3, 7, 2, 6]

[3, 7, 4, 5]

[3, 7, 4, 6]

[3, 7, 5, 4]

[3, 7, 5, 6]

[3, 7, 6, 4]

[3, 7, 6, 5]

[3, 7, 8, 5]

[3, 8, 2, 7]

[3, 8, 4, 5]

[4, 1, 2, 5]

[4, 1, 2, 6]

[4, 1, 3, 5]

[4, 1, 3, 6]

[4, 1, 5, 6]

[4, 1, 6, 5]

[4, 1, 7, 5]

[4, 2, 1, 5]

[4, 2, 1, 6]

[4, 2, 3, 5]

[4, 2, 3, 6]

[4, 2, 5, 6]

[4, 2, 6, 5]

[4, 3, 1, 5]

[4, 3, 1, 6]

[4, 3, 2, 5]

[4, 3, 2, 6]

[4, 3, 5, 6]

[4, 3, 6, 5]

[4, 3, 7, 6]

[4, 5, 1, 6]

[4, 5, 2, 6]

[4, 5, 3, 6]

[4, 5, 3, 7]

[4, 5, 6, 7]

[4, 5, 7, 6]

[4, 5, 8, 6]

[4, 6, 1, 5]

[4, 6, 2, 5]

[4, 6, 2, 7]

[4, 6, 3, 5]

[4, 6, 3, 7]

[4, 6, 5, 7]

[4, 6, 7, 5]

[4, 6, 7, 8]

[4, 6, 8, 5]

[4, 6, 8, 7]

[4, 7, 1, 5]

[4, 7, 2, 6]

[4, 7, 3, 5]

[4, 7, 3, 6]

[4, 7, 5, 6]

[4, 7, 6, 5]

[4, 7, 6, 8]

[4, 7, 8, 6]

[4, 8, 5, 7]

[4, 8, 6, 7]

[4, 8, 7, 5]

[4, 8, 7, 6]

[5, 1, 2, 6]

[5, 1, 3, 6]

[5, 1, 4, 6]

[5, 2, 1, 6]

[5, 2, 3, 6]

[5, 2, 4, 6]

[5, 2, 7, 6]

[5, 3, 1, 6]

[5, 3, 2, 6]

[5, 3, 4, 6]

[5, 3, 6, 7]

[5, 3, 7, 6]

[5, 4, 1, 6]

[5, 4, 2, 6]

[5, 4, 3, 6]

[5, 4, 6, 7]

[5, 4, 7, 6]

[5, 4, 8, 6]

[5, 6, 2, 7]

[5, 6, 3, 7]

[5, 6, 4, 7]

[5, 6, 8, 7]

[5, 7, 2, 6]

[5, 7, 3, 6]

[5, 7, 4, 6]

[5, 8, 3, 7]

[5, 8, 4, 6]

[5, 8, 6, 7]

[5, 8, 7, 6]

[6, 2, 5, 7]

[6, 3, 4, 7]

[6, 3, 5, 7]

[6, 4, 3, 7]

[6, 4, 5, 7]

[6, 4, 7, 8]

[6, 4, 8, 7]

[6, 5, 2, 7]

[6, 5, 3, 7]

[6, 5, 4, 7]

[6, 5, 8, 7]

[6, 7, 4, 8]

[6, 8, 2, 7]

[6, 8, 4, 7]

[7, 4, 6, 8]

[7, 6, 4, 8]

\end{multicols}


\begin{thebibliography}{56}

\bibitem{muniz} Mu{\~n}iz, A, Reverse the List of Numbers, Discussion on Mathstodon,  Available online:
\url{https://mathstodon.xyz/@two_star/112242224494626411} (accessed on 30/07/2026)

\bibitem{hn} Mu{\~n}iz, A, Reverse the List of Numbers, Discussion on Hacker News, Available online: \url{https://news.ycombinator.com/item?id=40010066} (accessed on 30/07/2026)

\bibitem{lawson-perfect} Puzzle implementation by Christian Lawson-Perfect, Available online: \url{https://somethingorotherwhatever.com/items/reverse-the-list-of-integers/} (accesssed on 31/07/2026)

\bibitem{alayachew} Puzzle implementation and solving algorithms by David Alayachew, Available online: \url{https://github.com/davidalayachew/ReverseTheListOfIntegers} (accessed on 31/07/2026)

\bibitem{61u} Puzzle solver by Hacker News user 61u, Available online: \url{https://ideone.com/MGFwdo} (accessed on 31/07/2026)

\bibitem{rokicki} Rokicki, T., Maximum number of moves to solve ``Reverse The List of Integers'' game with high value n, Entry A372008 in The On-line Encyclopedia of Integer Sequences, Available online: \url{https://oeis.org/A372008} (accessed on 30/07/2026)

\end{thebibliography}
\end{document}